\documentclass{gtart_h}  

\def\ifplaintex{\expandafter\ifx\csname documentclass\endcsname\relax}

\def\ifplaintex{\expandafter\ifx\csname documentclass\endcsname\relax}

\ifplaintex 
\else
\fi

\expandafter\ifx\csname epsfbox\endcsname\relax\input epsf\fi

\def\gt{{\mathsurround=0pt\it $\cal G\mskip-2mu$eometry \&\ 
$\cal T\!\!$opology}}        

\def\gtp{{\mathsurround=0pt\it $\cal G\mskip-2mu$eometry \&\ 
$\cal T\!\!$opology $\cal P\!$ublications}}  

\def\lognumber#1{\def\thelognumber{#1}}
\def\volumenumber#1{\def\thevolumenumber{#1}}
\def\papernumber#1{\def\thepapernumber{#1}}
\def\volumeyear#1{\def\thevolumeyear{#1}}

\def\pagenumbers#1#2{\def\startpage{#1}\def\finishpage{#2}}
\def\published#1{\def\publishdate{#1}}
\def\proposed#1{\def\theproposer{#1}}
\def\seconded#1{\def\theseconders{#1}}
\def\received#1{\def\receiveddate{#1}}

\def\accepted#1{\def\accepteddate{#1}}

\long\def\asciiabstract#1{\long\def\theasciiabstract{#1}}
\def\asciikeywords#1{\def\theasciikeywords{#1}}

\let\\\par\let\thelognumber\relax
\let\thevolumenumber\relax\let\thepapernumber\relax
\let\thevolumeyear\relax\let\thesamplenumber\relax\let\startpage\relax
\let\finishpage\relax\let\publishdate\relax\let\receiveddate\relax
\let\reviseddate\relax\let\accepteddate\relax\let\theasciititle\relax
\let\theasciiauthors\relax
\let\theasciiabstract\relax\let\theasciikeywords\relax
\let\theasciiemail\relax\let\theshortauthors\relax\let\theshorttitle\relax

\long\def\maketitlep{   

\count0=\startpage

\gt\hfill      
\hbox to 77pt{\vbox to 0pt{\vglue -15pt\epsfbox{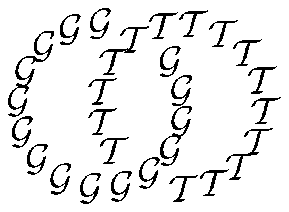}\vss}\hss}
\break
{\small\ifx\thesamplenumber\relax 
Volume \else Sample
\fi\thevolumenumber\ (\thevolumeyear)
\startpage--\finishpage\nl
Published: \publishdate}
\vglue 0.5truein plus 0.4fil minus 0.1truein

{\parskip=0pt\leftskip 0pt plus 1fil\def\\{\par\smallskip}{\ifplaintex\large
\else\Large\fi\bf\thetitle}\par\medskip}   

\vglue 0pt plus 0.1fil 

{\parskip=0pt\leftskip 0pt plus 1fil\def\\{\par}{\sc\theauthors}
\par\medskip}

\vglue 0pt plus 0.1fil 

{\small\parskip=0pt\let\newline\\
{\leftskip 0pt plus 1fil\def\\{\par}{\sl\theaddress}\par}
\expandafter\ifx\theemail\relax    
\relax\else\vglue 5pt plus 0.02fil minus 2pt\def\\{\stdspace{\rm 
and}\stdspace} 
\cl{Email:\stdspace\tt\theemail}\fi
\ifx\theurl\relax                  
\relax\else\vglue 5pt plus 0.02fil minus 2pt\def\\{\stdspace{\rm 
and}\stdspace}
\cl{URL:\stdspace\tt\theurl}\fi\par}

\vglue 7pt plus 0.3fil minus 3pt

{\bf Abstract}
\vglue 5pt plus 0.1fil minus 2pt

\theabstract

\vglue 7pt plus 0.3fil minus 3pt

{\bf AMS Classification numbers}\quad Primary:\quad \theprimaryclass

Secondary:\quad \thesecondaryclass

\vglue 5pt plus 0.3fil minus 2pt

{\bf Keywords:}\quad \thekeywords

\vglue 10pt plus 0.5fil minus 5pt

{\small  Proposed: \theproposer\hfill Received: \receiveddate\nl
Seconded: \theseconders\hfill 
\ifx\reviseddate\relax                         
Accepted: \accepteddate                        
\else
Revised: \reviseddate                          
\fi}
\eject
}       

\font\phead=cmsl9 scaled 950
\font=cmsl9 scaled 1050
\font\pnum=cmbx10 scaled 913
\font\lnum=cmbx10 
\font\pfoot=cmsl9 scaled 950
\font=cmsl9 scaled 1050
\ifplaintex
\headline{\vbox to 0pt{\vskip -4.5mm\line{\small\phead\ifnum
\count0=\startpage ISSN 1364-0380 (on line)
1465-3060 (printed) \hfill {\pnum\folio}\else\ifodd\count0\def\\{ }%
\ifx\theshorttitle\relax\thetitle\else\theshorttitle\fi\hfill{\pnum\folio}
\else\def\\{ and }{\pnum\folio}\hfill\ifx\theshortauthors\relax\theauthors
\else\theshortauthors\fi\fi\fi}\vss}}
\footline{\vbox to 0pt{\vglue 0mm\line{\small\pfoot\ifnum\count0=\startpage
\copyright\ \gtp\hfill\else
\gt, Volume \thevolumenumber\ (\thevolumeyear)\hfill\fi}\vss
}}
\else
\makeatletter
\def\@oddhead{{\small\lhead\ifnum\count0=\startpage ISSN 1364-0380 (on line)
1465-3060 (printed) \hfill {\lnum\number\count0}\else\ifodd\count0
\def\\{ }\ifx\theshorttitle\relax \thetitle \else\theshorttitle\fi\hfill
{\lnum\number\count0}\else\def\\{ and }{\lnum\number\count0}
\hfill\ifx\theshortauthors\relax 
\theauthors\else\theshortauthors\fi\fi\fi}}\def\@evenhead{\@oddhead}
\def\@oddfoot{\small\lfoot\ifnum\count0=\startpage\copyright\ \gtp\hfill\else
\gt, Volume \thevolumenumber\ (\thevolumeyear)\hfill\fi}
\def\@evenfoot{\@oddfoot}
\makeatother
\fi

\newwrite\gtoutfile
\long\gdef\makeheadfile{  
{\def\\{, }\def\s{ }
\immediate\openout\gtoutfile head.xxx
\immediate\write\gtoutfile{Proxy-for: \ifx\theasciiauthors\relax
\theauthors\else\theasciiauthors\fi\s<\ifx\theasciiemail\relax\theemail\else\theasciiemail\fi>}
\immediate\write\gtoutfile{\noexpand\\}
\immediate\write\gtoutfile{Authors: \ifx\theasciiauthors\relax
\theauthors\else\theasciiauthors\fi}
{\def\\{ }\immediate\write\gtoutfile{Title: \ifx\theasciititle\relax
\thetitle\else\theasciititle\fi}}
\immediate\write\gtoutfile{Subj-class: GT or SG or MG etc}
\immediate\write\gtoutfile{MSC-class: \theprimaryclass\ifx\thesecondaryclass\relax\else, \thesecondaryclass\fi}
\immediate\write\gtoutfile{Journal-ref: Geom. Topol. \thevolumenumber
(\thevolumeyear) \startpage-\finishpage}
\immediate\write\gtoutfile{Comments: Published by Geometry and Topology at}
\immediate\write\gtoutfile{\s\s http://www.maths.warwick.ac.uk/gt/GTVol\thevolumenumber/paper\thepapernumber.abs.html}
\immediate\write\gtoutfile{\noexpand\\}
\immediate\write\gtoutfile{}
\ifx\theasciiabstract\relax
\immediate\write\gtoutfile{\theabstract}\else
\immediate\write\gtoutfile{\theasciiabstract}\fi
\immediate\write\gtoutfile{}
\immediate\write\gtoutfile{\noexpand\\}
\immediate\write\gtoutfile{}
\immediate\closeout\gtoutfile}}  

\def\maketitlepage{\maketitlep\makeheadfile}
\let\maketitle\maketitlepage

\lognumber{537}
\received{7 December 2004}
\volumenumber{9}\papernumber{6}\volumeyear{2005}
\pagenumbers{203}{217}   
\published{20 January 2005}
\accepted{18 January 2005}
\proposed{Tomasz Mrowka}
\seconded{Ronald Fintushel, Ronald Stern}

\usepackage{amsmath,amssymb,pstricks}

\def\S{Section }

\newtheorem{theorem}{Theorem}
\newtheorem{lemma}[theorem]{Lemma}

\newtheorem{corollary}[theorem]{Corollary}
\theoremstyle{definition}
\newtheorem{definition}[theorem]{Definition}

\newcommand{\Z}{\mathbb{Z}}
\newcommand{\R}{\mathbb{R}}
\newcommand{\crit}{\mathrm{crit}}
\newcommand{\Map}{\mathrm{Map}{}}
\newcommand{\FA}{\mathcal{A}}
\newcommand{\FB}{\mathcal{B}}
\newcommand{\FC}{\mathcal{C}}
\newcommand{\FD}{\mathcal{D}}
\newcommand{\FF}{\mathcal{F}_0}

\begin{document}

\title{A stable classification of Lefschetz fibrations}
\author{Denis Auroux}

\address{Department of Mathematics, MIT\\Cambridge MA 02139, USA}
\email{auroux@math.mit.edu}

\begin{abstract}
We study the classification of Lefschetz fibrations up to stabilization
by fiber sum operations. We show that for each genus there is a
``universal'' fibration $f^0_g$ with the property that, if two Lefschetz
fibrations over $S^2$ have the same Euler--Poincar\'e characteristic and
signature, the same numbers of reducible singular fibers of each type,
and admit sections with the same self-intersection, then after repeatedly
fiber summing with $f^0_g$ they become isomorphic. As a consequence, 
any two compact integral symplectic 4--manifolds with the same values of
$(c_1^2,\,c_2,\,c_1\!\cdot\![\omega],\,[\omega]^2)$ become symplectomorphic
after blowups and symplectic sums with $f^0_g$.
\end{abstract}

\asciiabstract{%
We study the classification of Lefschetz fibrations up to
stabilization by fiber sum operations. We show that for each genus
there is a `universal' fibration f^0_g with the property that, if two
Lefschetz fibrations over S^2 have the same Euler-Poincare
characteristic and signature, the same numbers of reducible singular
fibers of each type, and admit sections with the same
self-intersection, then after repeatedly fiber summing with f^0_g they
become isomorphic. As a consequence, any two compact integral
symplectic 4-manifolds with the same values of (c_1^2, c_2, c_1.[w],
[w]^2) become symplectomorphic after blowups and symplectic sums with
f^0_g.}

\primaryclass{57R17}
\secondaryclass{53D35}
\keywords{Symplectic 4--manifolds, Lefschetz fibrations, fiber sums,
mapping class group factorizations}
\asciikeywords{Symplectic 4-manifolds, Lefschetz fibrations, fiber sums,
mapping class group factorizations}

\maketitle

\section{Introduction}

Lefschetz fibrations have been the focus of a lot of attention ever
since it was shown by Donaldson that, after blow-ups, every compact symplectic
4--manifold admits such structures \cite{Do}. We recall the definition:

\begin{definition}
A Lefschetz fibration on an oriented compact smooth 4--mani\-fold $M$ is a
smooth map $f\co M\to S^2$ which is a submersion everywhere except at finitely
many non-degenerate critical points $p_1,\dots,p_r$, near which $f$
identifies in local orientation-preserving complex coordinates with
the model map $(z_1,z_2)\mapsto z_1^2+z_2^2$.
\end{definition}

The smooth fibers of $f$ are compact surfaces, and the singular fibers
present nodal singularities; each singular fiber is obtained by collapsing
a simple closed loop (the {\it vanishing cycle}) in the smooth fiber.
The monodromy of the fibration around a singular fiber is given by a
positive Dehn twist along the vanishing cycle.

Denoting by $q_1,\dots,q_r\in S^2$ the images of the critical points (which
we will always assume to be distinct), and choosing a reference point 
$q_*\in S^2\setminus \crit(f)$, we can characterize the fibration $f$
by its {\it monodromy homomorphism}
$$\psi\co \pi_1(S^2\setminus\{q_1,\dots,q_r\},q_*)\to \Map_g,$$ where
$\Map_g=\pi_0\mathrm{Diff}^+(\Sigma_g)$ is the mapping class group of a
genus $g$ surface. It is a classical result (cf.\ \cite{Kas}) that the
monodromy morphism $\psi$ is uniquely determined up to conjugation by an
element of $\Map_g$ and the action of a braid on
$\pi_1(S^2\setminus\{q_i\})$ by ``Hurwitz moves'' (see \S \ref{s:fact});
moreover,
if the fiber genus is at least $2$ then the monodromy determines the
isomorphism class of the Lefschetz fibration $f$.

The classification of Lefschetz fibrations is a difficult problem
(essentially as difficult as the classification of symplectic 4--manifolds),
and is only understood in genus 1 and 2 (with some assumptions on the nature
of singular fibers in the latter case). It is a classical
result of Moishezon and Livne \cite{Mo1} that genus~1 Lefschetz fibrations
are all holomorphic, and are classified by the number of irreducible
singular fibers (which is a multiple of 12) and the number of reducible
singular fibers. More recently, Siebert and Tian \cite{ST} have obtained a
classification result for genus 2 Lefschetz fibrations without reducible
singular fibers and with ``transitive monodromy'' (a technical assumption
which we will not discuss here). Namely, these fibrations are all
holomorphic, and are classified by their
number of vanishing cycles, which is always a multiple of 10. In fact,
all such fibrations can be obtained as fiber sums of two standard
holomorphic fibrations $f_0$ and $f_1$ with respectively 20 and 30
singular fibers.

In higher genus, or even in genus 2 if one allows reducible singular fibers,
the classification appears to be much more complicated. However, we can
attempt to determine a minimal set of moves (ie, surgery operations) which
can be used to relate to each other any two Lefschetz fibrations with the
same genus. In this context, we consider stabilization by {\it fiber
sums} with certain standard fibrations. (The fiber sum of two Lefschetz
fibrations is obtained by deleting a neighborhood of a smooth fiber in each
of them, and gluing the resulting open manifolds along their boundaries in
a fiber-preserving manner). It was shown in \cite{Agen2} that,
given two genus 2 symplectic Lefschetz fibrations $f,f'$ with
the same numbers of singular fibers of each type (irreducible, reducible
with genus 1 components, reducible with components of genus 0 and 2),
for all large $n$ the fiber sums $f\#nf_0$ and $f'\#nf_0$ are isomorphic.
More generally,
as a corollary of a recent result of Kharlamov and Kulikov about braid
monodromy factorizations \cite{KK}, a similar result holds for all Lefschetz
fibrations with monodromy contained in the hyperelliptic mapping class
group.

Our goal is to obtain a similar stabilization result in the general case
(without assumptions on the fiber genus or on the monodromy).
In this context we must consider pairs of Lefschetz fibrations $f,f'$
with the same fiber genus and the same numbers of singular fibers of each
type (irreducible, or reducible of {\it type} $(h,g-h)$, ie,\ with
components of genera $h$ and $g-h$, for
each $0\le h\le \frac{g}{2}$), but we must also place two additional
restrictions (which automatically hold when $g\le 2$ or in the hyperelliptic
case). Namely, we must assume that the intersection forms on the total
spaces $M$ and $M'$ have the same signature, and we must assume that
the fibrations $f$ and $f'$ admit distinguished sections $s,s'$ which
represent classes in $H_2(M,\Z)$ (resp.\ $H_2(M',\Z)$) with
the same self-intersection number $-k$.

Then, we claim that, after repeatedly fiber summing $f$ and $f'$ with
a certain ``universal'' Lefschetz fibration $f^0_g$, constructed in
\S \ref{s:univ}, we
eventually obtain isomorphic Lefschetz fibrations:

\begin{theorem}\label{thm:main}
For every $g$ there exists a genus $g$ Lefschetz fibration $f^0_g$
with the following property.
Let $f\co M\to S^2$ and $f'\co M'\to S^2$ be two genus $g$ Lefschetz fibrations,
each equipped with a distinguished section. Assume that:

{\rm(i)}\qua the total spaces $M$ and $M'$ have the same Euler characteristic
and signature; 

{\rm(ii)}\qua the distinguished sections of $f$ and $f'$ have the
same self-intersection; 

{\rm(iii)}\qua $f$ and $f'$ have the same numbers of
reducible fibers of each type.

Then, for all large enough values of $n$, the fiber sums $f\# n\,f^0_g$ and
$f'\# n\,f^0_g$ are isomorphic.
\end{theorem}

A brief remark is in order about assumptions (i) and (ii) in this
statement. First, since the Euler characteristic is given by the formula
$\chi=4-4g+r$, where $g$ is the fiber genus and $r$ is the total number
of singular fibers, the first part of (i) is equivalent to
the requirement that $f$ and $f'$ have the same numbers of singular
fibers. Moreover, in the hyperelliptic case the assumption on signature
can be eliminated, because the signature is given by Endo's formula~\cite{E},
which involves only the number of singular fibers of each type; 
however, in general the signature depends on the actual vanishing cycles.
It is also worth mentioning that, in general, it is not known whether every
Lefschetz fibration admits a section (although there are no known examples
without a section). However, all Lefschetz fibrations obtained by blowing up
the base points of a pencil (and in particular all those which arise from
Donaldson's construction) admit sections of square $-1$.

The cases $g=0$ and $g=1$ of Theorem \ref{thm:main} are trivial
(in that case no stabilization is needed),
and the case $g=2$ is proved in \cite{Agen2} (taking $f^0_2$ to
be the holomorphic genus 2 fibration with 20 singular fibers and total
space a rational surface). Thus we will only consider the case $g\ge 3$
in the proof.

As a corollary of Theorem \ref{thm:main} and of Donaldson's result, we have
the following statement for {\em integral}\/ compact symplectic 4--manifolds
(ie, such that $[\omega]\in H^2(X,\R)$ is the image
of an integer cohomology class):

\begin{corollary}
Let $X,\,X'$ be two integral compact symplectic 4--manifolds with the same
$(c_1^2,\,c_2,\,c_1\mbox{$\cdot$}[\omega],\,[\omega]^2)$. Then $X$ and $X'$ become
symplectomorphic after sufficiently many blowups and symplectic sums with
the total space $X^0_g$ of the fibration $f^0_g$ (for a suitable genus $g$).
\end{corollary}

This corollary follows from Theorem \ref{thm:main} by considering pencils of
the same (large) degree $d$ on $X$ and $X'$, and blowing up the
$d^2[\omega]^2$ base points. The resulting
Lefschetz fibrations have the same fiber genus (by the assumptions on
$c_1\mbox{$\cdot$}[\omega]$ and $[\omega]^2$), admit sections of square
$-1$, and, if $d$ is large enough, can be assumed to contain only
irreducible fibers.

The proof of Theorem \ref{thm:main} actually gives a complete classification
of Lefschetz fibrations up to fiber sum stabilization. For example,
considering only Lefschetz fibrations with irreducible fibers, we have:

\begin{theorem}\label{thm:irredcase}
For every $g\ge 3$ there exist Lefschetz fibrations $f^A_g,f^B_g,f^C_g,f^D_g$
with the following property: if $f$ is a genus $g$ Lefschetz fibration
without reducible singular fibers, and if $f$ admits a section,
then there exist integers $a,b,c,d\in\Z$
such that for all large enough values of $n$ the fiber sums $f\#n\,f^0_g$
and $(n+a)\,f^A_g\,\#\,(n+b)\,f^B_g\,\#\,(n+c)\,f^C_g\,\#\,(n+d)\,f^D_g$ are
isomorphic.
\end{theorem}

The Lefschetz fibrations $f^A_g,f^B_g,f^C_g,f^D_g$ are constructed in
\S \ref{s:univ} (and $f^0_g$ is in fact nothing but their fiber sum).

The rest of this paper is organized as follows: in \S \ref{s:fact} we
review the description of Lefschetz fibrations by mapping class group
factorizations; in \S \ref{s:univ} we introduce the concept of universal
positive factorization and construct the Lefschetz fibrations $f^0_g$;
and in Sections \ref{s:stab}--\ref{s:proof} we prove Theorem \ref{thm:main}.
\medskip

{\sl This work was partially supported by NSF grant DMS-0244844.}

\section{Mapping class group factorizations}\label{s:fact}

The monodromy of a Lefschetz fibration can be encoded in a {\it mapping
class group factorization} by choosing an ordered system of generating loops
$\gamma_1,\dots,\gamma_r$ for $\pi_1(S^2\setminus\{q_1,\dots,q_r\})$, such that each
loop  $\gamma_i$ encircles only one of the points $q_i$ and $\prod \gamma_i$
is homotopically trivial. The monodromy of the fibration along each of
the loops $\gamma_i$ is a Dehn twist $\tau_i$; we can then describe the
fibration in terms of the relation $\tau_1\cdot \ldots\cdot \tau_r=1$ in
$\Map_g$. The choice of the loops $\gamma_i$ (and therefore of the twists
$\tau_i$) is of course not unique, but any two choices differ by a sequence
of {\it Hurwitz moves} exchanging consecutive factors:
$\tau_i\cdot \tau_{i+1} \to (\tau_{i+1})_{\tau_i^{-1}}\cdot
\tau_i$ or $\tau_i\cdot \tau_{i+1} \to \tau_{i+1}\cdot
(\tau_i)_{\tau_{i+1}}$, where we use the notation
$(\tau)_{\phi}=\phi^{-1}\tau \phi$, ie,\ if $\tau$ is a Dehn twist
along a loop $\delta$ then $(\tau)_\phi$ is the Dehn twist along the
loop $\phi(\delta)$.

\begin{definition}
A {\em factorization} $F=\tau_1\cdot\ldots\cdot \tau_r$ in $\Map_g$ is an
ordered tuple of positive Dehn twists. We say that two factorizations are 
{\em Hurwitz equivalent} $(F\sim F')$ if they can be obtained from each 
other by a sequence of Hurwitz moves.
\end{definition}

A Lefschetz fibration is thus characterized by a
factorization of the identity element in $\Map_g$, uniquely determined up 
to Hurwitz equivalence and simultaneous
conjugation of all factors by a same element of $\Map_g$, ie,\ up to
the equivalence relation generated by the moves
$$\tau_1\cdot\ldots\cdot\tau_i\cdot\tau_{i+1}\cdot\ldots\cdot\tau_r
\longleftrightarrow \tau_1\cdot\ldots\cdot\tau_{i+1}\cdot(\tau_i)_{\tau_{i+1}}
\cdot\ldots\cdot\tau_r\quad\forall 1\le i<r,$$
$$\tau_1\cdot\ldots\cdot\tau_r\longleftrightarrow (\tau_1)_\phi\cdot\ldots
\cdot(\tau_r)_\phi
\quad \forall \phi\in\Map_g.\smallskip$$

We will actually be considering Lefschetz fibrations equipped with a
distinguished section. The section determines a marked point in each fiber,
and so we can lift the monodromy to a {\it relative} mapping class group.
In fact, even though the normal bundle to the section $s$
is not trivial (it has degree $-k$ for some $k\ge 1$),
we can restrict ourselves to the
preimage of a large disc $\Delta$ containing all the critical values of $f$,
and fix a trivialization of the normal bundle to $s$ over $\Delta$.
Deleting a small tubular neighborhood of the section $s$, we can now view
the monodromy of $f$ as a homomorphism
$$\psi\co \pi_1(\Delta\setminus \{q_1,\dots,q_r\})\to \Map_{g,1},$$
where $\Map_{g,1}$ is the mapping class group of a genus $g$ surface with
one boundary component.
The product of the Dehn twists $\tau_i=\psi(\gamma_i)$ is not the identity,
but the central element $T_\delta^k\in\Map_{g,1}$, where $T_\delta$ is the
{\it boundary twist}, ie, the Dehn twist along a loop parallel to the
boundary.

With this understood, a Lefschetz fibration with a distinguished section
of square $-k$ is described by a factorization of $T_\delta^k$ as a product of
positive Dehn twists in $\Map_{g,1}$, up to Hurwitz equivalence and global
conjugation.
\medskip

A word about notations: while we use the multiplicative
notation for factorizations, and sometimes write $\tau_1\cdot\ldots\cdot
\tau_r=T_\delta^k$ to express the fact that $\tau_1\cdot\ldots\cdot\tau_r$
is a factorization of $T_\delta^k$, it is important not to confuse a
factorization (a tuple of Dehn twists) with the product of
its factors (an element in $\Map_{g,1}$). We will also use multiplicative
notation for the concatenation of factorizations ($F\cdot F'$ is the
factorization consisting of the factors in $F$, followed by those in $F'$,
and $(F)^n$ is the concatenation of $n$ copies of $F$),
and we will denote by $(F)_\phi$ the factorization obtained by conjugating
each factor of $F$ by the element $\phi\in\Map_{g,1}$.\medskip

To finish this section, we establish the following properties of Hurwitz
equivalence for factorizations of central elements:

\begin{lemma}\label{l:central}
Let $T$ be a central element in a group $G$. Then:

{\rm(a)}\qua if $F'\cdot F''$ is a factorization of $T$, then $F''\cdot
F'$ is also a factorization of $T$, and $F'\cdot F''\sim F''\cdot F'$;

{\rm(b)}\qua if $F$ is a factorization of $T$ whose factors generate $G$,
then $F\sim (F)_\phi$ for all $\phi\in G$;

{\rm(c)}\qua if $F$ is a factorization of $T$,
and $F'$ is any factorization, then $F'\cdot F\sim F\cdot F'$.
\end{lemma}

\proof (see also Lemma 6 in \cite{Agen2}).

(a)\qua To prove that any cyclic permutation of the factors amounts to a Hurwitz
equivalence, it suffices to prove that if $\tau\in G$ and
$\tau\cdot F''$ is a factorization
of $T$ then $\tau\cdot F''\sim F''\cdot \tau$. Denote by $\phi$ the product
of the factors in $F''$: using Hurwitz moves to move all the factors
in $F''$ to the left of $\tau$, we have $\tau\cdot F''\sim F''\cdot
(\tau)_\phi$. The result then follows from the observation that
$\phi=\tau^{-1}T$ commutes with $\tau$.

(b)\qua Let $\tau$ be any of the factors in $F$: then by (a) we can perform a
cyclic permutation of the factors and obtain a
factorization $F'$ such that $F\sim F'\cdot \tau$. Moving
$\tau$ to the left of $F'$, we have $F'\cdot \tau\sim \tau\cdot
(F')_\tau=(\tau\cdot F')_\tau$. Applying (a) again we have $(\tau\cdot
F')_\tau\sim (F)_\tau$. So, for any factor $\tau$ of $F$, we have
$F\sim (F)_\tau$, and similarly $F\sim (F)_{\tau^{-1}}$. The result then
follows from the assumption that the factors of $F$ generate $G$, by
expressing $\phi$ in terms of the factors.

(c)\qua Simply move all the factors of $F$ to the left of the factors in $F'$, to
obtain $F'\cdot F\sim F\cdot (F')_T=F\cdot F'$ (since $T$ is central).
\endproof

\section{Universal positive factorizations}\label{s:univ}

Let us first recall a presentation of $\Map_{g,1}$ due to Matsumoto
\cite{Ma}, which is a reformulation of Wajnryb's classical presentation
\cite{Wa} in a form that is more convenient for our purposes
(see Theorem 1.3 and Remark 1.1 of \cite{Ma}):

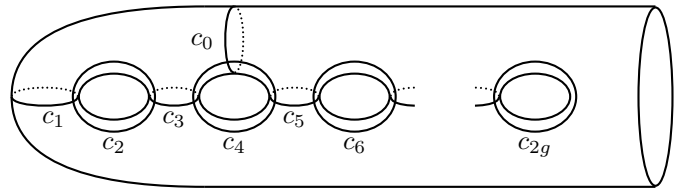
\begin{figure}[b]\small
\centerline{
\setlength{\unitlength}{8mm}
\begin{picture}(11,3)(-0.5,-1.5)
\psset{unit=\unitlength}
\psbezier{-}(-0.2,0)(-0.2,1.25)(1.5,1.5)(3,1.5)
\psbezier{-}(-0.2,0)(-0.2,-1.25)(1.5,-1.5)(3,-1.5)
\psline{-}(3,1.5)(10.5,1.5)
\psline{-}(3,-1.5)(10.5,-1.5)
\psellipse(10.5,0)(0.3,1.5)
\psellipse(1.5,0)(0.6,0.4)
\psellipse(3.5,0)(0.6,0.4)
\psellipse(5.5,0)(0.6,0.4)
\psellipse(8.5,0)(0.6,0.4)
\psbezier{-}(-0.2,0)(-0.2,-0.2)(0.9,-0.2)(0.9,0)
\psbezier[linestyle=dotted,dotsep=1pt]{-}(-0.2,0)(-0.2,0.2)(0.9,0.2)(0.9,0)
\psbezier{-}(2.1,0)(2.1,-0.2)(2.9,-0.2)(2.9,0)
\psbezier[linestyle=dotted,dotsep=1pt]{-}(2.1,0)(2.1,0.2)(2.9,0.2)(2.9,0)
\psbezier{-}(4.1,0)(4.1,-0.2)(4.9,-0.2)(4.9,0)
\psbezier[linestyle=dotted,dotsep=1pt]{-}(4.1,0)(4.1,0.2)(4.9,0.2)(4.9,0)
\psbezier{-}(6.1,0)(6.1,-0.15)(6.5,-0.15)(6.5,-0.15)
\psbezier[linestyle=dotted,dotsep=1pt]{-}(6.1,0)(6.1,0.15)(6.5,0.15)(6.5,0.15)
\psbezier{-}(7.9,0)(7.9,-0.15)(7.5,-0.15)(7.5,-0.15)
\psbezier[linestyle=dotted,dotsep=1pt]{-}(7.9,0)(7.9,0.15)(7.5,0.15)(7.5,0.15)
\psellipse(1.5,0)(0.7,0.6)
\psellipse(3.5,0)(0.7,0.6)
\psellipse(5.5,0)(0.7,0.6)
\psellipse(8.5,0)(0.7,0.6)
\psbezier{-}(3.5,0.4)(3.3,0.4)(3.3,1.5)(3.5,1.5)
\psbezier[linestyle=dotted,dotsep=1pt]{-}(3.5,0.4)(3.7,0.4)(3.7,1.5)(3.5,1.5)
\put(0.5,-0.25){\makebox(0,0)[ct]{$c_1$}}
\put(2.5,-0.25){\makebox(0,0)[ct]{$c_3$}}
\put(4.5,-0.25){\makebox(0,0)[ct]{$c_5$}}
\put(1.5,-0.68){\makebox(0,0)[ct]{$c_2$}}
\put(3.5,-0.68){\makebox(0,0)[ct]{$c_4$}}
\put(5.5,-0.68){\makebox(0,0)[ct]{$c_6$}}
\put(8.5,-0.68){\makebox(0,0)[ct]{$c_{2g}$}}
\put(3.2,0.8){\makebox(0,0)[rb]{$c_0$}}
\end{picture}
\qquad
}
\caption{The Dehn--Lickorish--Humphries generators of $\Map_{g,1}$}
\label{fig:mapg1}
\end{figure}

\begin{theorem}[Matsumoto]\label{thm:present}
For $g\ge 2$, the mapping class group $\Map_{g,1}$ is generated by the
Dehn twists $a_0,\dots,a_{2g}$ along the loops $c_0,\dots,c_{2g}$
represented in Figure \ref{fig:mapg1}, with the relations:

{\rm(i)}\qua $a_ia_j=a_ja_i$ if $c_i\cap c_j=\emptyset$, and
$a_ia_ja_i=a_ja_ia_j$ if $c_i\cap c_j\neq\emptyset$;

{\rm(ii)}\qua $(a_0a_2a_3a_4)^{10}=(a_0a_1a_2a_3a_4)^6$;

{\rm(iii)}\qua for $g\ge 3:$
$(a_0a_1a_2a_3a_4a_5a_6)^{9}=(a_0a_2a_3a_4a_5a_6)^{12}$.
\end{theorem}

The relations (i) are the {\it braid relations}, and
realize $\Map_{g,1}$ as a quotient of an Artin group, while (ii)
is a reformulation of the {\it chain relation}, and (iii) is a
reformulation of the {\it lantern relation} (see \cite{Ma}).

The subgroup of $\Map_{g,1}$ generated by $a_1,\dots,a_{2g}$ is the
{\it hyperelliptic} subgroup, and is closely related to the braid group
$B_{2g+1}$ (realizing the genus $g$ surface as a double cover of the disc
branched in $2g+1$ points, the Dehn twists $a_1,\dots,a_{2g}$ are the lifts
of the standard generators of $B_{2g+1}$).

\begin{lemma}\label{l:rn}
For every integer $1<n<2g$, let
\begin{multline*}
R_n=(a_{n+1}\cdot a_{n+2}\cdot\ldots\cdot a_{2g})^{2g-n+1}\cdot\\
\cdot \prod_{i=n}^1
(a_i\cdot a_{i+1}\cdot\ldots\cdot a_{i+2g-n})\cdot
\prod_{i=2g-n+1}^1 (a_i\cdot a_{i+1}\cdot\ldots\cdot a_{i+n-1}).
\end{multline*}
Then
$(a_1\cdot\ldots\cdot a_{n-1})^{2n}\cdot (R_n)^2$ is a factorization of
$T_\delta$ in $\Map_{g,1}$.
\end{lemma}

\proof
We work in the braid group $B_{2g+1}$ with generators $x_1,\dots,x_{2g}$,
and consider the expression obtained from $R_n$ after replacing
each $a_i$ by $x_i$. Then it is easy to see that $a'=(x_1\dots x_{n-1})^n$
is the full twist rotating the $n$ leftmost strands by $2\pi$, while
$a''=(x_{n+1}\dots x_{2g})^{2g-n+1}$ is the full twist rotating the
$2g+1-n$ rightmost strands by $2\pi$. Moreover, $b'=\prod_{i=n}^1
(x_i\dots x_{i+2g-n})$ is the braid which exchanges the $n$ leftmost strands
with the $2g-n+1$ rightmost strands in the counterclockwise direction, while
$b''=\prod_{i=2g-n+1}^1 (x_i\dots x_{i+n-1})$ does the same for the $2g-n+1$
leftmost strands and the $n$ rightmost strands. Hence the product $b'b''$
corresponds to a full rotation of the $n$ leftmost strands around the 
$2g-n+1$ rightmost strands, and $a'a''b'b''$ is the
full twist $\Delta^2=(x_1\dots x_{2g})^{2g+1}$.
Since $\Delta^2$ is a central element in $B_{2g+1}$, we also have
$a''b'b''a'=\Delta^2$.

We now lift things to the double cover; since $\Delta^2$ lifts
to the hyperelliptic element $H$ (rotating the surface about its central
axis by $\pi$), we deduce from the above calculation that
$(a_1\cdot\ldots\cdot a_{n-1})^n\cdot R_n$ and $R_n\cdot (a_1\cdot\ldots
\cdot a_{n-1})^n$ are factorizations of $H$, and hence that $(a_1\cdot\ldots
\cdot a_{n-1})^n\cdot (R_n)^2\cdot (a_1\cdot\ldots\cdot a_{n-1})^n$ is a
factorization of $H^2=T_\delta$. Since $T_\delta$ is central in
$\Map_{g,1}$, the result follows by Lemma \ref{l:central}(a).
\endproof

It is in fact not hard to check explicitly that the factorization considered
in Lemma \ref{l:rn} is Hurwitz equivalent to the standard hyperelliptic
factorization \hbox{$(a_1\cdot\ldots\cdot a_{2g})^{4g+2}$}.
\medskip

From now on we assume that $g\ge 3$. By Theorem 1.4
of \cite{Ma}, $(a_0a_2a_3a_4)^{10}=(a_0a_1a_2a_3a_4)^6
=(a_1a_2a_3a_4)^{10}$ and $(a_0a_1a_2a_3a_4a_5a_6)^9=
(a_0a_2a_3a_4a_5a_6)^{12}=(a_1a_2a_3a_4a_5a_6)^{14}$.
Hence, we can define new factorizations of
$T_\delta$ by substitution into the factorization of Lemma \ref{l:rn}:

\begin{definition}
Let $\FA=(a_0\cdot a_2\cdot a_3\cdot a_4)^{10}\cdot (R_5)^2$, 
$\FB=(a_0\cdot a_1\cdot a_2\cdot a_3\cdot a_4)^{6}\cdot (R_5)^2$,
$\FC=(a_0\cdot a_1\cdot a_2\cdot a_3\cdot a_4\cdot a_5\cdot a_6)^{9}\cdot
(R_7)^2$,
$\FD=(a_0\cdot a_2\cdot a_3\cdot a_4\cdot a_5\cdot a_6)^{12}\cdot (R_7)^2$
(where for $g=3$ we take $R_7$ to be the empty factorization),
and $\FF=\FA\cdot\FB\cdot\FC\cdot\FD.$
\end{definition}

$\FA,\,\FB,\,\FC,\,\FD$ are factorizations of the central element $T_\delta$
in which every factor is one of the $(a_i)_{0\le i\le 2g}$, and every
generator appears at least once (except possibly for $\FD$, which does
not involve $a_1$ when $g=3$). 

We also define $f^0_g,\,f^A_g,\,f^B_g,\,f^C_g,\,f^D_g$ to be the Lefschetz
fibrations with monodromy factorizations $\FF,\,\FA,\,\FB,\,\FC,\,\FD$
respectively (so $f^A_g,\,f^B_g,\,f^C_g,\,f^D_g$ are irreducible and admit
sections of square $-1$, while $f^0_g$ is their fiber sum and admits a section of
square $-4$). Let us mention that, as a consequence of Lemma
\ref{l:central}(b), when performing a fiber sum with $f^0_g$ the choice of
the identification diffeomorphism between fibers is irrelevant, and all
possible ways in which the fiber sum can be carried out are equivalent.

The factorizations $\FA,\FB,\FC,\FD$ form a ``universal'' set of positive
factorizations, in
the sense that their factors are exactly the generators of $\Map_{g,1}$
(out of sequence, and with some repetitions), and every relation in
the presentation of Theorem \ref{thm:present} can be interpreted either as
a Hurwitz equivalence or as a substitution replacing one of these
factorizations by another one of them. We will see below
that these properties are the key ingredients for the proof of Theorem
\ref{thm:main}; since many other groups related to braid groups or mapping
class groups can be presented in a similar manner, the methods
used here may also be relevant to the study of factorizations in these groups.

\section{Stable equivalence of factorizations}\label{s:stab}

In this section, we prove the following result, which implies Theorem
\ref{thm:irredcase}:

\begin{theorem}\label{thm:stab}
Let $F,F'$ be two factorizations of the same element of\/ $\Map_{g,1}$
as a product of positive Dehn twists along non-separating curves. Then
there exist integers $a,b,c,d,k,l$ such that $F\cdot (\FA)^a\cdot (\FB)^b
\cdot (\FC)^c\cdot (\FD)^d\sim F'\cdot (\FA)^{a+l}\cdot (\FB)^{b-l}\cdot
(\FC)^{c+k}\cdot(\FD)^{d-k}.$
\end{theorem}

In order to prove this result, we consider factorizations where the factors
are either positive Dehn twists or their inverses, and the
equivalence relation $\equiv$ generated by the following moves:

$\bullet$ Hurwitz moves involving only positive Dehn twists;

$\bullet$ creation or cancellation of pairs of inverse factors:
$a_i\cdot a_i^{-1}\equiv a_i^{-1}\cdot a_i\equiv \emptyset$;

$\bullet$ defining relations of the mapping class group: $a_i\cdot a_j
\equiv a_j\cdot a_i$ if $c_i\cap c_j=\emptyset$, $a_i\cdot a_j\cdot a_i
\equiv a_j\cdot a_i\cdot a_j$ if $c_i\cap c_j\neq\emptyset$, 
$(a_0\cdot a_2\cdot a_3\cdot a_4)^{10}\equiv (a_0\cdot a_1\cdot a_2\cdot
a_3\cdot a_4)^{6}$, and $(a_0\cdot a_1\cdot a_2\cdot a_3\cdot a_4\cdot
a_5\cdot a_6)^9\equiv (a_0\cdot a_2\cdot a_3\cdot a_4\cdot a_5\cdot
a_6)^{12}.$

\begin{lemma}\label{l:simplify}
If the factors of $F$ are Dehn twists along non-separating curves,
then there exists a factorization $\bar{F}$ in which every factor is of
the form $a_i^{\pm 1}$ for some $0\le i\le 2g$, and such that $F\equiv
\bar{F}$.
\end{lemma}

\proof
We use pair creations and Hurwitz moves to replace
every factor in $F$ by a factorization involving only the $a_i^{\pm 1}.$
Let $\tau$ be a factor in $F$. Since $\tau$ is a Dehn twist along a
non-separating
curve, there exist $g_0,\dots,g_k\in \{a_0^{\pm 1},\dots,a_{2g}^{\pm 1}\}$
such that $\tau=(\prod_1^k g_j)^{-1}g_0\,(\prod_1^k g_j)$. We proceed by
induction on $k$. If $k=0$ then $\tau$ is already one of the generators.
Otherwise, if $g_k$ is one of the generators, say $a_i$, then we can write
$\tau=(a_i^{-1}\tilde\tau a_i)\equiv a_i^{-1}\cdot a_i\cdot (a_i^{-1}\tilde
\tau a_i)\equiv a_i^{-1} \cdot \tilde\tau\cdot a_i$ (using a pair creation
and a Hurwitz move). Similarly, if $g_k=a_i^{-1}$, then we can write
$\tau=(a_i\tilde\tau a_i^{-1})\equiv (a_i\tilde\tau a_i^{-1})\cdot a_i\cdot
a_i^{-1}\equiv a_i\cdot \tilde\tau\cdot a_i^{-1}.$ Since $\tilde\tau$ is
conjugated to one of the generators by a word of length $k-1$, this
completes the proof.
\endproof

\begin{lemma}\label{l:equiv}
Under the assumptions of Theorem \ref{thm:stab}, $F\equiv F'$.
\end{lemma}

\proof
We first use Lemma \ref{l:simplify} to replace $F$ and $F'$ by equivalent
factorizations $\bar{F}$ and $\bar{F}'$ whose factors are all of the form
$a_i^{\pm 1}$. Next, recall that if a group $G$ admits a presentation with
generators $\{a_i,\ i\in I\}$ and relations $\{r_j,\ j\in J\}$, then it is
generated as a {\it monoid}\/ by the elements $\{a_i,a_i^{-1},\ i\in I\}$,
and a presentation of $G$ as a monoid is given by the set of relations
$R'=\hbox{$\{r_j,\ j\in J\}$}\cup \{a_ia_i^{-1}=1,\,a_i^{-1}a_i=1,\ i\in I\}$. Hence,
if $\bar{F}$ and $\bar{F}'$ are factorizations of the same element,
then we can rewrite one into the other by successively applying the
rewriting rules given by the set of relations $R'$. However,
in the case of the
mapping class group, each rewriting is one of the moves that generate
the equivalence relation $\equiv$ (either one of the defining
relations of $\Map_{g,1}$, or the creation or cancellation of a pair of
inverses).
\endproof

Denote by $\equiv^+$ the equivalence relation generated by Hurwitz moves
and by the defining relations, ie,\ without allowing creations of pairs of
inverse factors. Then we have:

\begin{lemma}\label{l:stequiv}
Under the assumptions of Theorem \ref{thm:stab}, there exists an integer $n$
such that $F\cdot (\FA)^n\equiv^+ F'\cdot (\FA)^n$.
\end{lemma}

\proof
By Lemma \ref{l:equiv}, $F\equiv F'$, so we can transform $F$ into $F'$ by
a sequence of Hurwitz moves, pair creations/cancellations, and defining
relations. Call $F=F_0,F_1,\dots,F_m=F'$ the successive factorizations
appearing in this sequence of moves; let $n_j$ be the number of
factors of the form $a_i^{-1}$ appearing in the factorization $F_j$, and
let $n=\sup\{n_0,\dots,n_m\}$.

Recall that the factors of $\FA$ generate $\Map_{g,1}$; therefore, by Lemma
\ref{l:central}(a), for every $i\in\{0,\dots,2g\}$ there exists a
factorization $\FA_i$ whose factors are elements of $\{a_0,\dots,a_{2g}\}$,
and such that $\FA\sim a_i\cdot \FA_i\sim \FA_i\cdot a_i$. (For example
$\FA_i$ can be obtained by cyclically permuting the factors of $\FA$ and
deleting an occurrence of $a_i$). 
Let $\smash{F_j^+}$ be the factorization
obtained from $F_j$ by replacing each factor of the form $a_i^{-1}$ by the
factorization $\FA_i$. Then we claim that, for all $0\le j<m$, 
$F_j^+\cdot (\FA)^{n-n_j}\equiv^+ F_{j+1}^+\cdot (\FA)^{n-n_{j+1}}.$

Indeed, if $F_{j+1}$ is obtained from $F_j$ by a Hurwitz move or by applying
a defining relation, then the negative factors are not involved and the
claim is obvious. If $F_{j+1}$ is obtained from $F_j$ by deleting a pair
of mutually inverse factors $a_i\cdot a_i^{-1}$, $F_{j+1}^+$ is obtained
from $F_j^+$ by deleting an occurrence of the subword 
$a_i\cdot \FA_i$. Hence, we can write $F_j^+=F'_j\cdot a_i\cdot \FA_i\cdot
F''_j$ and $F_{j+1}^+=F'_j\cdot F''_j$ for some $F'_j,F''_j$,
and the claim follows from the sequence of Hurwitz moves
$$F'_j\cdot a_i\cdot \FA_i\cdot F''_j\cdot
(\FA)^{n-n_j}\sim F'_j\cdot \FA\cdot F''_j\cdot (\FA)^{n-n_j}\sim
F'_j\cdot F''_j\cdot(\FA)^{n-n_j+1},$$
where in the last step we have used Lemma \ref{l:central}(c).
The argument is the same for creations of pairs of inverses.
The proof is then completed by observing that $F_0^+=F$ and $F_m^+=F'$,
since $F$ and $F'$ contain no negative factors.
\endproof

We can now proceed with the proof of Theorem \ref{thm:stab}.
By Lemma \ref{l:stequiv}, there exists $n$ such that 
$F\cdot (\FA)^n\equiv^+ F'\cdot (\FA)^n$, so we can transform
$F\cdot(\FA)^n$ into $F'\cdot(\FA)^n$ by a sequence of Hurwitz moves
and applications of the defining relations.
Let $F_0=F\cdot (\FA)^n,\,F_1,\,\dots,\,F_m=F'\cdot(\FA)^n$ be the
successive factorizations appearing in this sequence of moves. If $F_{j+1}$
is obtained from $F_j$ by a Hurwitz move, or by applying one of
the braid relations, then we have
$F_{j+1}\sim F_j$. For example, a braid relation
of the form $a_i\cdot a_j\cdot a_i\equiv a_j\cdot a_i\cdot a_j$ can be
viewed as a succession of two Hurwitz moves $a_i\cdot a_j\cdot a_i\sim
a_j\cdot (a_i)_{a_j}\cdot a_i\sim a_j\cdot a_i\cdot (a_i)_{a_ja_i}$, where
$(a_i)_{a_ja_i}=\smash{(a_ja_i)^{-1}}a_ia_ja_i=a_j$. 

On the other hand, if
$F_{j+1}$ is obtained from $F_j$ by applying the relation (ii) from Theorem
\ref{thm:present}, then we can write $F_j=F'_j\cdot(a_0\cdot a_2\cdot a_3\cdot
a_4)^{10}\cdot F''_j$ for some $F'_j,F''_j$, and $F_{j+1}=F'_j
\cdot(a_0\cdot a_1\cdot a_2\cdot a_3\cdot a_4)^{6}\cdot F''_j$.
It is then easy to check, using Lemma \ref{l:central}\,(a) and (c), that
$F_j\cdot \FB\sim F_{j+1}\cdot \FA$; and vice-versa if we apply relation
(ii) backwards. Similarly, if $F_{j+1}$ is obtained from $F_j$ by applying
relation (iii), then $F_j\cdot \FD\sim F_{j+1}\cdot \FC$, and vice-versa
if we apply relation (iii) backwards.

Hence, if we concatenate each $F_j$ with suitable numbers of copies
of $\FA$, $\FB$, $\FC$ and $\FD$ (depending on $j$), then we can realize
each step as a Hurwitz equivalence. Since we always trade a copy of
$\FA$ for a copy of $\FB$, and a copy of $\FC$ for a copy of $\FD$,
Theorem \ref{thm:stab} follows.
\medskip

We can now prove Theorem \ref{thm:irredcase}:

\proof[Proof of Theorem \ref{thm:irredcase}]
Let $F$ be a factorization in $\Map_{g,1}$ associated to the Lefschetz
fibration $f$: then the product of the factors in $F$ is equal to
$T_\delta^m$, for some integer $m\ge 1$ (such that the chosen section of $f$
has self-intersection $-m$). The result then follows
by applying Theorem \ref{thm:stab} to $F$ and $F'=(\FA)^m$.
\endproof

\section{Proof of Theorem \ref{thm:main}}\label{s:proof}

Let $F$ and $F'$ be factorizations in $\Map_{g,1}$ describing the
monodromies of the Lefschetz fibrations $f$ and $f'$. Assumption (ii)
on the self-intersection numbers of the distinguished sections implies
that the products of the factors in $F$ and $F'$ are equal to each other,
and are of the form $T_\delta^m$ for some $m\ge 1$.
We first deal with the reducible singular fibers, using the following lemma:

\begin{lemma}\label{l:red}
If $\tau,\tau'$ are Dehn twists along separating curves of the same type,
then there exists an integer $n$ and a factorization $F''$ involving only
Dehn twists along non-separating curves, such that $\tau\cdot(\FA)^n\sim
\tau'\cdot F''$.
\end{lemma}

\proof
$\tau,\tau'$ are conjugated to each other in $\Map_{g,1}$, so there exist
$g_1,\dots,g_k\in\{a_0^{\pm 1},\dots,a_{2g}^{\pm 1}\}$ such that
$\tau'=(\prod_1^k g_j)^{-1}\tau\,(\prod_1^k g_j).$ It is enough to consider
the case $k=1$ (iterating $k$ times in the general case).
If $\tau'=a_i^{-1}\tau a_i$ then, with the same
notations as in the proof of Lemma \ref{l:stequiv}, we have $\tau\cdot \FA
\sim \tau\cdot a_i\cdot \FA_i\sim a_i\cdot \tau'\cdot \FA_i\sim \tau'\cdot
(a_i)_{\tau'}\cdot \FA_i$, and the result follows by setting
$F''=(a_i)_{\tau'}\cdot \FA_i$. Similarly, if $\tau'=a_i\tau a_i^{-1}$ then
$\tau\cdot \FA\sim \FA\cdot \tau\sim \FA_i\cdot a_i\cdot \tau\sim \FA_i
\cdot \tau'\cdot a_i\sim \tau'\cdot (\FA_i)_{\tau'}\cdot a_i.$
\endproof

The manner in which we use this lemma is the following: let $s$ be the
number of reducible singular fibers of $f$ and $f'$. Without loss of
generality, we can assume that the homologically trivial vanishing cycles
correspond to the first $s$ factors of $F$ and $F'$, and that they are
ordered according to types (this can always be ensured by performing Hurwitz
moves). Call these factors $\tau_1,\dots,\tau_s$ for $F$, and
$\tau'_1,\dots,\tau'_s$ for $F'$. Then assumption (iii) on the numbers of
reducible singular fibers implies that $\tau_j$ and $\tau'_j$ are conjugated
for each $1\le j\le s$. Hence, applying Lemma \ref{l:red} to each pair
$(\tau_j,\tau'_j)$, and adding sufficiently many copies of $\FA$ to $F$
(using Lemma \ref{l:central}(c) to move them to the beginning of the
factorization), we can replace each $\tau_j$ by $\tau'_j$, at the expense
of generating extra Dehn twists along nonseparating curves. After suitable
Hurwitz moves, we conclude that there exists an integer $N$ and
factorizations $\tilde{F},\tilde{F}'$ involving only Dehn twists along
non-separating curves, such that $F\cdot (\FA)^N\sim
\tau'_1\cdot\ldots\cdot\tau'_s\cdot \tilde{F}$
and $F'\cdot (\FA)^N\sim \tau'_1\cdot\ldots\cdot\tau'_s\cdot\tilde{F}'.$

Since $\tilde{F}$ and $\tilde{F'}$ are factorizations of the same element
$(\tau'_1\dots\tau'_s)^{-1}T_\delta^{m+N}$, we can apply Theorem
\ref{thm:stab} to them. It follows that there exist integers $a,b,c,d,k,l$
such that $F\cdot (\FA)^{N+a}\cdot (\FB)^b\cdot (\FC)^c\cdot (\FD)^d
\sim F'\cdot (\FA)^{N+a+l}\cdot (\FB)^{b-l}\cdot (\FC)^{c+k}\cdot
(\FD)^{d-k}$. This implies that the fiber sums
$\hat{f}=f\,\#\,(N+a)f^A_g\,\#\,bf^B_g\,\#\,cf^C_g\,\#\,df^D_g$
and $\hat{f}'=f'\,\#\,(N+a+l)f^A_g\,\#\,(b-l)f^B_g\,\#\,(c+k)f^C_g\,\#\,(d-k)f^D_g$ are isomorphic.
Performing additional fiber sums if necessary, we can assume
that $N+a=b=c=d$. Then, in order to complete the proof of
Theorem~\ref{thm:main}, it is sufficient to prove that $k=l=0$.
For this purpose we use the following lemmas to compare the Euler--Poincar\'e
characteristics and signatures of the total spaces $\hat{M}$ and $\hat{M}'$
of $\hat{f}$ and $\hat{f}'$:

\begin{lemma}
$\chi(\hat{M}')-\chi(\hat{M})=\chi(M')-\chi(M)+10\,l-9\,k$.
\end{lemma}

\proof Recall that the Euler characteristic of a genus $g$ Lefschetz
fibration over $S^2$ with $r$ singular fibers is equal to $4-4g+r$. Hence,
we just have to compare the numbers of singular fibers of $\hat{f}$ and
$\hat{f}'$. Since $f^A_g$ has 10 more singular fibers than $f^B_g$, and
$f^C_g$ has 9 fewer singular fibers than $f^D_g$, the result follows.
\endproof

\begin{lemma}
$\sigma(\hat{M}')-\sigma(\hat{M})=\sigma(M')-\sigma(M)-6\,l+5\,k$.
\end{lemma}

\proof By Novikov additivity, it is sufficient to show that the signatures
of the total spaces $M_A,M_B,M_C,M_D$ of $f^A_g,f^B_g,f^C_g,f^D_g$ satisfy
the relations $\sigma(M_A)=\sigma(M_B)-6$ and $\sigma(M_C)=\sigma(M_D)+5$.

These signatures can be computed explicitly via an algorithm due to Ozbagci
\cite{Oz}. Since Ozbagci's formula is a sum of individual contributions
which each depend only on one of the factors and on the {\it product} of all the
preceding factors, it is sufficient to carry out the algorithm for the
portions of $\FA$ and $\FB$ (resp.\ $\FC$ and $\FD$) which differ from each
other; the contributions from the common part $(R_5)^2$ (resp.\ $(R_7)^2$)
will be the same in both cases.

In fact, after a closer look at the signature formula it is easy to convince
oneself that
$\sigma(M_A)-\sigma(M_B)$ and $\sigma(M_C)-\sigma(M_D)$ do not depend on
$g$, and can be computed for a fixed low value of $g$ (e.g., $g=3$).

Then, rather than Ozbagci's somewhat complicated formula, one can use
the following simple recipe to determine the signature -- the underlying
principle being that, given a Lefschetz fibration $f\co M\to S^2$ admitting
a section, the complement to the fiber and section classes in
$H_2(M,\Z)$ is generated by certain linear combinations of the Lefschetz
thimbles of $f$.

Given the set of vanishing cycles $(\delta_1,\dots,\delta_r)$ (ie,
loops in the fiber $\Sigma_g$ such that each monodromy factor $\tau_i$
is the Dehn
twist along $\delta_i$), form the $r\times r$ matrix $Q$ whose entries
are given by $$q_{ij}=\begin{cases} 0 & \mbox{if }i>j,\\
-1 & \mbox{if }i=j,\\
\delta_i\cdot\delta_j &\mbox{if }i>j,\end{cases}$$
where $\delta_i\cdot\delta_j$ is the intersection number in
$H_1(\Sigma_g,\Z)$. In a suitable sense, $Q$ is the matrix of the
intersection pairing on the space of formal linear combinations of Lefschetz
thimbles, and its antisymmetrization $A=Q-Q^t$ describes the intersection
pairing between vanishing cycles inside $\Sigma_g$.

Viewing $Q$ and $A$ as bilinear forms, the kernel of $A$ is the space of
all combinations of Lefschetz thimbles which have homologically trivial
boundary in $H_1(\Sigma_g,\Z)$, and can hence be completed to 2--cycles inside
$M$. The restriction $Q'=Q_{|\mathrm{Ker}\,A}$ is now a (degenerate)
symmetric bilinear form, of rank $b_2(M)-2$; and $Q'$ has the same
signature as the intersection form on $H_2(M,\Z)$, ie\ $\sigma(Q')=\sigma(M)$.

Applying this formula, we easily check that for $g=3$, $\sigma(M_A)=-48$,
$\sigma(M_B)=-42$, $\sigma(M_C)=-35$, and $\sigma(M_D)=-40$.
\endproof

The proof of Theorem \ref{thm:main} can now be completed by observing that,
since $\chi(M')=\chi(M)$ and $\sigma(M')=\sigma(M)$ by assumption (i),
and since $\hat{M}$ and $\hat{M}'$ are diffeomorphic by construction, we must have
$10\,l=9\,k$ and $6\,l=5\,k$, which implies that $k=l=0$.

\end{document}